\documentstyle{amsart}

\newtheorem{theorem}{Theorem}[section]
\newtheorem{defn}[theorem]{Definition}
\newtheorem{prop}[theorem]{Proposition}
\newtheorem{lem}[theorem]{Lemma}

\newtheorem{claim}[theorem]{Claim}
\newcommand{\sub}{\subseteq}

\newcommand{\ov}{\overline}
\newcommand{\os}{\overline{s}}

\newcommand{\rt}{\rightarrow}
\newcommand{\sm}{\setminus}
\newcommand{\proof}{\noi {\bf Proof:}\ }
\newcommand{\noi}{\noindent}
\newsymbol\upharpoonright 1316
\newcommand{\rest}{\upharpoonright}
\newcommand{\con}{\;^\wedge }
\newcommand{\ra}{\rangle}
\newcommand{\la}{\langle}

\newcommand{\CF}{{\cal F}}
\newcommand{\CG}{{\cal G}}
\newcommand{\CO}{{\cal O}}
  
\newcommand{\CT}{{\cal T}}   
\newcommand{\CX}{{\cal X}} 
\newcommand{\CY}{{\cal Y}}  
\newcommand{\CZ}{{\cal Z}}  
\newcommand{\BI}{{\bf I }}
\newcommand{\BII}{{\bf II }}
\newcommand{\Fr}{{\frak Fr }}
\newcommand{\om}{\omega}

\newcommand{\ffo}{\; ^{<\omega} \! \omega}
\newcommand{\fso}{[\omega]^{<\omega}}
\newcommand{\ff}{\; ^{\omega} \! 2}
\newcommand{\fff}{\; ^{<\omega} \! 2}
\newcommand{\ffso}{\; ^{<\omega} \!( [\omega]^{<\omega})}
\newcommand{\so}{[\omega]^\omega}
\newcommand{\po}{\mbox{{\Large $\wp$}}(\omega)}
\newcommand{\fa}{\forall}					  
\newcommand{\fai}{\forall^\infty}					  
\newcommand{\ex}{\exists}

\begin{document}

\title{Filter Games and Combinatorial Properties of Winning Strategies}

\author{Claude Laflamme}
\address{Department of Mathematics and Statistics \\
         University of Calgary \\
         Calgary, Alberta  \\
         Canada T2N 1N4}
\email{laflamme@@acs.ucalgary.ca}


\date{April 12, 1994}

\thanks{This research was partially supported by  NSERC of Canada.}

 \maketitle

\begin{abstract}

We characterize winning strategies in various infinite games involving
filters on the natural numbers in terms of combinatorics or structural
properties of the given filter.  These generalize several ultrafilter
games of Galvin. 

 \end{abstract}
							  
\section{Introduction}

We look at various infinite games between two players \BI and \BII
involving filters on the natural numbers in which \BI either plays
cofinite sets, members of $\CF$ or $\CF^+$ and player \BII responds with
an element or a finite subset of \BI's move, depending on the game.  In
each version, the outcome depends on the set produced by player \BII,
whether it belongs to the given filter $\CF$, $\CF^+$ or even $\CF^*$,
the dual ideal. 

In each game considered, we will characterize winning strategies of
either player in terms of combinatorics of the given filter $\CF$; these
combinatorics turn out to be generalizations of the classical notions of
P-points, Q-points and selectivity for ultrafilters.  In the case of
ultrafilters, $\CF = \CF^+$ and most of our games to ones already
studied by Galvin in unpublished manuscripts \cite{GA}; the various
generalized combinatorics enjoyed by the filters become equivalent. 

Several characterizations of Ramsey ultrafilters and P-points were known
from works of Booth (\cite{BO}) and Kunen (\cite{KU}), and some
generalizations of these combinatorics to filters were already made by
Grigorieff in \cite{GR} where for example the notion of P-filter is
characterized in terms of branches through certain trees; we shall see
that this is very much in the spirit of winning strategies for certain
games. 

Other variations of these games for ultrafilters can be found in
Chapter VI of Shelah's book \cite{S}, and two of the games below have
been analyzed by Bartoszynski and Scheepers \cite{BS}. 

We wish to thank Chris Leary for helpful suggestions and corrections
regarding the present paper.

\medskip

Our terminology is standard but we review the main concepts and
notation.  The natural numbers will be denoted by $\omega$, $\po$
denotes the collection of all its subsets.  Given $X \in \po$, we write
$[X]^\omega$ and $[X]^{<\omega}$ to denote the infinite or finite
subsets of $X$ respectively.  We use the well known `almost inclusion'
ordering between members of $\so$, i.e.  $X \sub^*Y$ if $X \sm Y$ is
finite.  We identify $\po$ with $\ff$ via characteristic functions.  The
space $\ff$ is further equipped with the product topology of the discrete
space $\{0,1\}$.  A basic neighbourhood is then given by sets of the
form

\[ \CO_s = \{f \in \ff: s \sub f \} \] 

\noi where $s \in \fff$, the collection of finite binary sequences.  The
terms ``nowhere dense'', ``meager'', ``Baire property'' all refer to
this topology. Concatenation of elements ${\ov s}, {\ov t} \in \ffo$ will
be written ${\ov s} \;^\wedge {\ov t}$.

A filter is a collection of subsets of $\omega$ closed under finite
intersections, supersets and containing all cofinite sets; it is called
proper if it contains only infinite sets.  For a filter $\CF$, $\CF^+$
denotes the collection of all sets $X$ such that $\la \CF,X \ra$ is a
proper filter; it is useful to notice that $X \in \CF^+$ if and only if
$X^c \notin \CF$.  $(\CF^+)^c = \wp(\omega ) \sm \CF^+$, the collection
of sets incompatible with $\CF$ is the dual ideal and is usually denoted
by $\CF^*$.  The {\em Fr\'echet} filter is the collection of cofinite
sets, denoted by $\Fr$. 

The families $\CF$ and $\CF^+$ are dual in a different sense; this means
that a set $X$ containing an element of each member of $\CF$ (resp. 
$\CF^+$) must belong to $\CF^+$ (resp.  $\CF$).  In particular $\Fr$ and
$\so$ are dual.  From more general work of Aczel (\cite{A}) and Blass
(\cite{B}), there is a duality between games in which a player chooses
$X_k \in \CF$ while the other player responds with $n_k \in X_k$, and
games in which a player chooses $Y_k \in \CF^+$ while the other player
responds with $n_k \in Y_k$. The point is that the statements

\[ (\fa X \in \CF)(\ex n \in X) \phi(n) \mbox{ and }
   (\ex Y \in \CF^+)(\fa n \in Y) \phi(n) \]

\noi are equivalent.

The following important result characterizes meager filters in terms
of combinatorial properties.

\begin{prop}\label{TAL}
(Talagrand (\cite{T})) The following are equivalent for a filter $\CF$:
\begin{enumerate}
\item $\CF$ has the Baire property.
\item $\CF$ is meager.
\item There is a sequence $n_0 < n_1 < \cdots$ such that
   \[ (\fa X \in \CF)(\fai k) \; X \cap [n_k,n_{k+1}) \neq \emptyset. \]
\end{enumerate}
\end{prop}

Combinatorial properties of filters have played an important role in
applications of Set Theory, and the classical notions of a filter being
meager, a P-filter or selective have been around a long time.  These
concepts will be generalized below in terms of trees and other structural
properties; these combinatorial ideas have their roots in Ramsey theory
and P-points and selective ultrafilters (sometimes called `Ramsey') have
characterizations in term of these trees; this can be found in the
papers by Booth \cite{BO} and Grigorieff \cite{GR}. 

We call a tree $\CT \sub \ffo$ an $\CX$-tree for some $\CX \sub \so$
($\CX$ will usually be a filter $\CF$ or $\CF^+$), if for each ${\ov s}
\in \CT$, there is an $X_{\ov s} \in \CX$ such that ${\ov s} \con n \in
\CT$ for all $n \in X_{\ov s}$.  Similarly we call a tree $\CT \sub \ffso$ an
$\CX$-tree of finite sets for some $\CX \sub \so$, if for each ${\ov s}
\in \CT$, there is an $X_{\ov s} \in \CX$ such that ${\ov s} \con a \in
\CT$ for each $a \in [X_{\ov s}]^{<\om}$.  A branch of such a tree is thus an
infinite sequence of finite sets and we will be interested in the union
of such a branch, an infinite subset of $\om$. 

Here are a few more combinatorial properties of filters that we will consider.

\begin{defn}
Let $\CF$  be a filter on $\omega$.
\begin{enumerate}

\item $\CF$ is called a Q-filter if for any partition of $\omega$ into finite
sets $\la s_k : k \in \omega \ra$, there is an $X \in \CF$ such that
$\mid X \cap s_k \mid \leq 1$ for all $k$.

\item $\CF$ is called a weak Q-filter if for any partition of $\omega$
into finite sets $\la s_k : k \in \omega \ra$, there is an $X \in \CF^+$
such that $\mid X \cap s_k \mid \leq 1$ for all $k$. 

\item $\CF$ is called {\it diagonalizable} if there is an $X \in \so$
such that $X \sub^* Y$ for all $Y \in \CF$. 

\item $\CF$ is called {\it $\omega$-diagonalizable} if there are $\la
X_n \in \so:n \in \omega \ra$ such that for each $Y \in \CF$, there is
an $n$ such that $X_n \sub^* Y$. 

\item $\CF$ is called {\it $\omega$-+-diagonalizable} if there are $\la
X_n \in \CF^+:n \in \omega \ra$ such that for each $Y \in \CF$, there is
an $n$ such that $X_n \sub^* Y$. 

\item A set $X \sub \fso$ is called $\CZ$-universal ($\CZ$ will be $\CF$
or $\CF^+$) if for each $Y \in \CZ$, there is an $x \in X \cap
[Y]^{<\omega}$.  $\CF$ is called {\em $\omega$-diagonalizable by
$\CZ$-universal sets} if there are $\CZ$-universal sets $\la X_n: n \in
\omega \ra$ such that for all $Y \in \CF$, there is an $n$ such that $x
\cap Y \neq \emptyset$ for all but finitely many $x \in X_n$. 

\item $\CF$ is a {\em P-filter} if given any sequence $\la X_n :n \in
\omega \ra \sub \CF$,there is an $X \in \CF$ such that $X \sub^* X_n$
for each $n$. 

\item $\CF$ is a {\em weak P-filter} if given any sequence $\la X_n:n
\in \omega \ra \sub \CF$, there is an $X \in \CF^+$ such that $X \sub^*
X_n$ for each $n$.  Equivalently, every $\CF$-tree of finite sets has a
branch whose union is in $\CF^+$.

\item $\CF$ is a {\em P$^+$-filter} if every $\CF^+$-tree of finite sets
has a branch whose union is if $\CF^+$. 

\item $\CF$ is {\em Ramsey} if any $\CF$-tree has a branch in $\CF$;
equivalently, $\CF$ is both a Q-filter and a P-filter. 

\item $\CF$ is {\em weakly Ramsey} if any $\CF$-tree has a branch
in $\CF^+$.

\item $\CF$ is {\em +-Ramsey} if every $\CF^+$-tree has a branch
in $\CF^+$.
 
\item $\CF$ is a P-point if it is an ultrafilter that is also
a P-filter.

\end{enumerate} 
\end{defn}

If $\CF$ is a P-filter, then $\CF$ is diagonalizable if and only if it
is $\omega$-diagonalizable if and only if it is
$\omega-+$-diagonalizable.  On the other hand, these notions are
distinct.  Indeed the filter 

\[\CF = \Fr \otimes \Fr = \{ X \sub \om \times \om :
\{n : \{m : (n,m) \in X \} \mbox{ is cofinite } \}   \mbox{ is cofinite } \} \]

\noi is $\om$-diagonalizable, but not $\om$-+-diagonalizable; similarly,
if $\CG$ is any non-dia\-go\-na\-li\-zable filter (any non-meager filter
will do), then the filter

\noi $ \CF = \{X \sub \om \times \om :
             \mbox{ for each }n, \{m : (n,m) \in X \}  \in \CG,$

\hfill and $ \{n : \{m : (n,m) \in X \} \mbox{ is cofinite } \} 
        \mbox{ is cofinite } \} $

\noi is $\om$-+-diagonalizable, but not diagonalizable.

If $\CF$ is an ultrafilter, then $\CF$ is a P-filter if and only if it
is a weak P-filter if and only if it is a P$^+$-filter; it is a Q-filter
if and only if it is a weak Q-filter, and $\CF$ is Ramsey if and only if it
is weakly Ramsey if and only if it is +-Ramsey.  Therefore these notions
generalize the classical combinatorial properties of ultrafilters, but
again, these notions can be seen to be different for filters  in general. 

The notions of $\omega$-diagonalizability by $\CF$-(resp
$\CF^+$)-universal sets are generalizations of the regular
$\omega$-(resp.  +)-diagonalizability.  Observe that a filter $\CF$ is
$\om$-diagonalizable by $\Fr$-universal sets if and only if
diagonalizable by a single $\Fr$-universal set if and only if it is
meager, and $\om$-+-diagonalizability implies $\om$-diagonalizability by
$\CF$-universal sets. This notion appears to be new.

Tree combinatorics is what most interest us in this paper as they
naturally occur in terms of winning strategies; our main effort is
then to express these combinatorics in terms of more familiar
concepts. The  following lemma shows the spirit of the paper.

\begin{lem}\label{NMP}
$\CF$ is a non-meager P-filter if and only if every $\CF$-tree of
finite sets has a branch whose union is in $\CF$.
\end{lem}

\proof Assume first that every $\CF$-tree of finite sets has a branch
whose union is in $\CF$. Given a descending sequence $\la A_n ; n \in
\omega \ra \sub \CF$, define an $\CF$-tree $\CT \sub \ffso$ such that
$X_{\overline s}=A_n$ for each ${\overline s} \in \CT \cap \; ^n \!(
[\omega]^{<\omega})$. Any branch through $\CT$ whose union is in $\CF$
shows that $\CF$ is a P-filter.  \noi To verify that $\CF$ is
non-meager, consider an increasing sequence of natural numbers $n_0 <
n_1 < \ldots$ and build again an $\CF$-tree of finite sets $\CT$ as
follows. Having already ${\overline s} \in \CT$, choose $k$ such that
${\overline s} \in \; ^{<\omega} \!( [n_k]^{<\omega})$ and let
$X_{\overline s} = \omega \setminus n_{k+1}$. The union of any branch
through $\CT$ misses infinitely many intervals of the form
$[n_k,n_{k+1})$ and therefore shows that $\CF$ is also non-meager.

Now assume that $\CF$ is a non-meager P-filter and let $\CT$ be an
$\CF$-tree of finite sets. Define $n_0=0$ and $A_0 =
 X_\emptyset$. More generally, given $A_0 \supseteq A_1 \supseteq
\cdots \supseteq A_k$ and $n_0 < n_1 < \cdots < n_k$, let $n_{k+1} >
n_k$ in $A_k$ and put

\[ A_{k+1} = \bigcap \{X_{\la s_0,s_1,\cdots ,s_i \ra }:
         i \leq k+1, s_j \sub A_j \cap [0,n_{k+1}] \}, \]

\noi Now as $\CF$ is a P-filter, there is a $Y \in \CF$ such that $Y \sub^*
A_n$ for each $n$ and we may as well assume by a reindexing that
\[ Y \sm n_{k+1} \sub A_k \mbox{ for each } k. \]

\noi Now as $\CF$ is also non-meager, we can find 
an infinite set $K=\{k_\ell : \ell \in \om \}$ and we might as well
assume that for each $\ell$,
 \[ Y \cap [ n_{k_\ell}, n_{k_\ell + 1}) = \emptyset . \]
\noi Define $s_\ell =  Y \cap [ n_{k_\ell}, n_{k_{\ell + 1}}) = 
              Y \cap [ n_{k_\ell + 1 }, n_{k_{\ell + 1}})$ for each $\ell$.
We claim that $\la s_k : k \in \omega \ra$ is a branch through $\CT$;
indeed 

\[\begin{array}{ll}
s_\ell   &  =  Y \cap [ n_{k_\ell}, n_{k_{\ell + 1}}) \\
         &  =  Y \cap [ n_{k_\ell + 1 }, n_{k_{\ell + 1}}) \\
      & \sub A_{k_\ell } \sub X_{\la s_0,s_1,\cdots, s_{\ell -1} \ra }.
\end{array}\]

\noi But this concludes the proof as its union $Y$ is in $\CF$.
\qed

\bigskip

\section{Filter Games}

We will be interested in infinite games of the form ${\frak
G}(\CX,\CY,\CZ)$ where $\CX$ will usually be a filter $\CF$ or $\CF ^+$,
$\CY$ will be either $\omega$ or $\fso$, and $\CZ$ will be either $\CF$,
$\CF^+$, $\CF^c$, the complement of $\CF$ in $\po$ or else $(\CF^+)^c=\CF^*$,
the dual ideal.

The game ${\frak G}(\CX,\CY,\CZ)$ is played by two players \BI and \BII
as follows: at stage $k<\omega$, \BI chooses $X_k \in \CX$, then \BII
responds with either $n_k \in X_k$ in the case that $\CY$ is $\omega$,
or else responds with a nonempty $s_k \in [X_k]^{<\omega}$ in the case
that $\CY$ is $\fso$.  At the end of the game, \BII is declared the
winner if $\{n_k:k \in \omega \} (\mbox{resp.  } \bigcup_{k \in \omega}
s_k) \in \CZ$. 
 
A few variations of these games have been considered in the literature,
see for example in \cite{BS} and \cite{S} for some special cases.  In
particular, the game ${\frak G}(\Fr, \omega,\CZ)$ is equivalent to the
game in which at stage $k$ player \BI chooses $m_k \in \omega$ and \BII
responds with $n_k > m_k$, the outcome being that \BII wins the play if
$\{n_k : k \in \omega \} \in \CZ$ as before.  Therefore we start with a
result of \cite{BS}. 

\begin{theorem}(\cite{BS})\label{FR-W-F} Fix a filter $\CF$ and consider
the game ${\frak G}(\Fr, \omega, \CF)$.  Then

\begin{enumerate}
  \item \BI has a winning strategy if and only if $\CF$ is not a Q-filter.
  \item \BII has no winning strategy.
\end{enumerate}

\end{theorem}

But this game generalizes in many ways and we have the following results.

\begin{theorem}\label{FR-W-F+} Fix a filter $\CF$ and consider
the game ${\frak G}(\Fr, \omega, \CF^+)$.  Then

\begin{enumerate}
  \item \BI has a winning strategy if and only if $\CF$ is not a weak Q-filter .
  \item \BII has a winning strategy 
             if and only if $\CF$ is $\omega$-diagonalized.
\end{enumerate}

\end{theorem}

\noi A much less interesting game is the following:

\begin{theorem}\label{FR-W-FC} Fix a filter $\CF$ and consider
the game ${\frak G}(\Fr, \omega, \CF^c)$.  Then the game is determined and

\begin{enumerate}
  \item  \BII always has a winning strategy.
\end{enumerate}

\end{theorem}

\begin{theorem}\label{FR-W-F+C} Fix a filter $\CF$ and consider
the game ${\frak G}(\Fr, \omega, \CF^*)$.  Then the game is
determined and

\begin{enumerate}
  \item  \BI has a winning strategy if and only if $\CF=\Fr$
          if and only if  \BII has no winning strategy.
\end{enumerate}

\end{theorem}

Because of the duality between $\Fr$ and $\Fr^+ = \so$, we will get

\begin{theorem}\label{DUAL-FR} Fix a filter $\CF$, then the following games
are dual of each other; that is a player has a winning strategy in one game
if and only if the other player has a winning strategy in the other game.

\begin{enumerate}
  \item ${\frak G}(\so, \omega, \CF)$ and ${\frak G}(\Fr, \omega, \CF^c)$ .
  \item ${\frak G}(\so, \omega, \CF^+)$ and ${\frak G}(\Fr, \omega, \CF^*)$.
  \item ${\frak G}(\so, \omega, \CF^c)$ and ${\frak G}(\Fr, \omega, \CF)$.
  \item ${\frak G}(\so, \omega, \CF^*)$ and ${\frak G}(\Fr, \omega, \CF^+)$.
\end{enumerate}

\end{theorem}

Now we consider the more interesting games where \BI plays members of
the filter $\CF$. 

\begin{theorem}\label{F-W-F} Fix a filter $\CF$ and consider
the game ${\frak G}(\CF, \omega, \CF)$.  Then

\begin{enumerate}
  \item \BI has a winning strategy if and only if $\CF$ is not a Ramsey filter.
  \item \BII never has a winning strategy.
\end{enumerate}

\end{theorem}

\begin{theorem}\label{F-W-F+} Fix a filter $\CF$ and consider
the game ${\frak G}(\CF, \omega, \CF^+)$.  Then

\begin{enumerate}
  \item \BI has a winning strategy if and only if
          $\CF$ is not  weakly Ramsey.
  \item \BII has a winning strategy if and only if
          $\CF$ is $\omega$-+-diagonalizable.
\end{enumerate}

\end{theorem}

\begin{theorem}\label{F-W-FC} Fix a filter $\CF$ and consider
the game ${\frak G}(\CF, \omega, \CF^c)$.  Then

\begin{enumerate}
  \item \BI never has a winning strategy.
  \item \BII has a winning strategy if and only 
         if $\CF$ is not a Ramsey ultrafilter.
\end{enumerate}

\end{theorem}

\begin{theorem}\label{F-W-F+C} Fix a filter $\CF$ and consider
the game ${\frak G}(\CF, \omega, \CF^*)$.  Then

\begin{enumerate}
  \item \BI has a winning strategy if and only if $\CF$ is countably generated.
  \item \BII has a winning strategy 
               if and only if $\CF$ is not a +-Ramsey filter.
\end{enumerate}

\end{theorem}

Again the duality of $\CF$ and $\CF^+$ will provide the following result.

\begin{theorem}\label{DUAL-F} Fix a filter $\CF$, then the following games
are dual of each other; that is a player has a winning strategy in one game
if and only if the other player has a winning strategy in the other game.

\begin{enumerate}
  \item ${\frak G}(\CF^+, \omega, \CF)$ and ${\frak G}(\CF, \omega, \CF^c)$ .
  \item ${\frak G}(\CF^+, \omega, \CF^+)$ and ${\frak G}(\CF, \omega, \CF^*)$.
  \item ${\frak G}(\CF^+, \omega, \CF^c)$ and ${\frak G}(\CF, \omega, \CF)$.
  \item ${\frak G}(\CF^+, \omega, \CF^*)$ and ${\frak G}(\CF, \omega, \CF^+)$.
\end{enumerate}

\end{theorem}

Before we turn to games in which player \BII chooses finite sets at each round,
consider the following infinite game ${\frak G}_1(\CF)$ defined in \cite{BS}:
at stage $k$, player \BI chooses $m_k \in \omega$ and \BII responds with $n_k$.
At the end, \BII is declared the winner if 
\begin{enumerate}
 \item $n_1 < n_2 < \cdots < n_k < \cdots$,
 \item $m_k < n_k$ for infinitely many $k$, and
 \item $\{n_k : k \in \omega \} \in \CF$. 
\end{enumerate} 

\noi It is proved in \cite{BS} that \BII does not have any winning
strategy in ${\frak G}_1(\CF)$ and that \BI has a winning strategy if
and only if $\CF$ is meager.  We have the following.

\begin{theorem}\label{FR-FSW-F} Fix a filter $\CF$, the the games
${\frak G}_1(\CF)$ and ${\frak G}(\Fr, \fso, \CF)$ are equivalent; that
is a player has a winning strategy in one game if and only the same
player has a winning strategy in the other game.  Therefore, by
\cite{BS}, we have for either game,

\begin{enumerate}
  \item \BI has a winning strategy if and only if $\CF$ is meager.
  \item \BII never has a winning strategy.
\end{enumerate}

\end{theorem}

\begin{theorem}\label{FR-FSW-F+} Fix a filter $\CF$, then the game
${\frak G}(\Fr, \fso, \CF^+)$ is dual to the game ${\frak G}(\Fr, \fso,
\CF)$.  Therefore

\begin{enumerate}
  \item \BI never has a winning strategy.
  \item \BII has a winning strategy if and only if $\CF$ is meager.
\end{enumerate}

\end{theorem}

As above, the following game is uninteresting.

\begin{theorem}\label{FR-FSW-FC} Fix a filter $\CF$, then the game
${\frak G}(\Fr, \fso, \CF^c)$ is equivalent to ${\frak G}(\Fr, \omega,
\CF^c)$.  Therefore

\begin{enumerate}
  \item  \BII always has a winning strategy.
\end{enumerate}

\end{theorem}

\begin{theorem}\label{FR-FSW-F+C} Fix a filter $\CF$, the 
the game ${\frak G}(\Fr, \fso , \CF^*)$ is equivalent to
 ${\frak G}(\Fr, \omega , \CF^*)$. Therefore

\begin{enumerate}
  \item  \BI has a winning strategy if and only if $\CF=\Fr$
          if and only if  \BII has no winning strategy.
\end{enumerate}

\end{theorem}

We turn to those games where \BI plays members of $\CF$ and \BII responds 
with finite subsets.

\begin{theorem}\label{F-FSW-F} Fix a filter $\CF$ and consider
the game ${\frak G}(\CF, \fso, \CF)$.  Then

\begin{enumerate}
  \item \BI has \underline{no} winning strategy if and only if $\CF$ is
          a non-meager P-filter.
  \item \BII never has a winning strategy.
\end{enumerate}

\end{theorem}

\begin{theorem}\label{F-FSW-F+} Fix a filter $\CF$ and consider
the game ${\frak G}(\CF, \fso, \CF^+)$.  Then

\begin{enumerate}
  \item \BI has \underline{no} winning strategy if and only if $\CF$ is a weak
          P-filter.
  \item \BII has a winning strategy 
     if and only if $\CF$ is $\omega$-diagonalizable  by $\CF$-universal sets.
\end{enumerate}

\end{theorem}

\begin{theorem}\label{F-FSW-FC} Fix a filter $\CF$, then
the game ${\frak G}(\CF, \fso, \CF^c)$ is equivalent to the game
 ${\frak G}(\CF, \omega, \CF^c)$. Therefore

\begin{enumerate}
  \item \BI never has a winning strategy.
  \item \BII has a winning strategy if and only if $\CF$ is not a Ramsey
             ultrafilter.
 \end{enumerate}

\end{theorem}

\begin{theorem}\label{F-FSW-F+C} Fix a filter $\CF$, then the 
game ${\frak G}(\CF, \fso, \CF^*)$ is equivalent to the game 
${\frak G}(\CF, \omega, \CF^*)$.  Thus

\begin{enumerate}
  \item \BI has a winning strategy if and only if $\CF$ is countably 
           generated.
  \item \BII has a winning strategy if and only if $\CF$ is not a Ramsey 
           ultrafilter.
 \end{enumerate}

\end{theorem}

Finally, we turn to games where \BI plays members of $\CF^+$ while 
\BII replies with finite subsets. Note that we do not have here the
same duality as when \BII responded with natural numbers.

\begin{theorem}\label{F+-FSW-F} Fix a filter $\CF$ and consider
the game ${\frak G}(\CF^+, \fso, \CF)$.  Then

\begin{enumerate}
  \item \BI has a winning strategy if and only if $\CF$ is not a P-point.
  \item \BII never has a winning strategy.
 \end{enumerate}

\end{theorem}

\begin{theorem}\label{F+-FSW-F+} Fix a filter $\CF$ and consider
the game ${\frak G}(\CF^+, \fso, \CF^+)$. 

\begin{enumerate}
  \item \BI has a winning strategy if and only if $\CF$ is not a P$^+$-filter.
  \item \BII has a winning strategy if and only if $\CF$ is $\omega$-diagonalizable
        by $\CF^+$-universal sets.
 \end{enumerate}

\end{theorem}

\begin{theorem}\label{F+-FSW-FC} Fix a filter $\CF$, then
the game ${\frak G}(\CF^+, \fso, \CF^c)$ is dual to
the game ${\frak G}(\CF, \omega, \CF)$. Therefore

\begin{enumerate}
  \item \BI never has a winning strategy.
  \item \BII has a winning strategy if and only if 
         $\CF$ is not a Ramsey filter.
 \end{enumerate}

\end{theorem}

\begin{theorem}\label{F+-FSW-F+C} Fix a filter $\CF$, then 
the game ${\frak G}(\CF^+, \fso, \CF^*)$ is dual to
the game ${\frak G}(\CF, \omega, \CF^+)$. Therefore

\begin{enumerate}
  \item \BI has a winning strategy if and only if
       $\CF$ is $\omega$-+-diagonalizable.
  \item \BII has a winning strategy if and only if
         $\CF$ is not weakly Ramsey.
 \end{enumerate}

\end{theorem}

\section{Proofs}

In this section we verify the results of section 2. We start with two
general results.

\begin{lem} \label{EQUI} 

If a family $\CZ \sub \po$ is closed under supersets, then the two games
${\frak G}(\CX, \omega, \CZ^c)$ and ${\frak G}(\CX, \fso, \CZ^c)$ are
equivalent. 

\end{lem}

\proof Since player \BII is trying to get out of $\CZ$ which is assumed
to be closed under supersets (therefore $\CZ^c$ is closed under
subsets), the best strategy for \BII is to play finite sets as small as
possible, namely singleton since a legal move must be nonempty.  \qed

The next Lemma regards the duality mentioned in the introduction and is
taken from the work of Aczel \cite{A} and Blass (see \cite{B}, Theorem
1). We include a hint of the proof for completeness.

\begin{theorem} \label{BLASS}

(\cite{A},\cite{B}) For a given filter $\CF$, the game ${\frak G}(\CF,
\omega, \CZ )$ and the game ${\frak G}(\CF^+, \omega, \CZ^c)$ are dual;
that is a player has a winning strategy in one of these games if and
only if the other player has a winning strategy in the other game. 

\end{theorem}

\proof Suppose \BII has a winning strategy \$ in the game 
${\frak G}(\CF, \omega, \CZ )$, we define a strategy \$\$ for \BI in 
the  game  ${\frak G}(\CF^+, \omega, \CZ^c)$ as follows:

\noi \BI starts with \$\$($\emptyset$)= $\{ \$(X): X \in \CF \} \in
\CF^+$.  When \BII responds with $n_0$, \BI remembers one set $X_0 \in
\CF$ such that \$$(X_0)=n_0$. 

\noi At stage $k$, \BI has remembered $k$ sets $X_0,X_1,\cdots,X_{k-1}$
from $\CF$ while \BII responded with $\la n_0,n_1,\cdots,n_{k-1} \ra$. 
\BI then plays 

\[\$\$( \la n_0,n_1,\cdots,n_{k-1} \ra ) =
   \{ \$( \la X_0,X_1,\cdots,X_{k-1},X \ra): X \in \CF \} \in \CF^+; \]

\noi \BII responds with $n_k$ and \BI remembers one set $X_k \in \CF$
such that 

\[ \$( \la X_0,X_1,\cdots,X_{k-1},X_k \ra )= n_k. \]

\noi Thus a play in the new game corresponds to a play in the former game and
thus the outcome $\{n_k : k \in \omega \} \in \CZ$ and \BI 's strategy is 
a winning strategy in  ${\frak G}(\CF^+, \omega, \CZ^c)$. 

The other cases are quite similar and left to the reader. \qed

Now we are ready to attack the proofs of section 2.

\bigskip

\noi \underline{{PROOF OF THEOREM  \ref{FR-W-F+}:}} We first deal
with player \BI.  So suppose that $\CF$ is not a weak Q-filter  and therefore
there is a partition of $\omega$ into finite sets $\la s_k:k\in \omega
\ra$ such that no $X \in \CF^+$ meet each $s_k$ in at most one point. 
Then \BI's strategy at stage $k$, after \BII has played $\la
n_0,n_1,\cdots,n_{k-1} \ra$, is to respond with $\bigcup \{s_i : s_i
\cap \{n_0,n_1,\cdots,n_{k-1}\} = \emptyset \}$. 

\noi Now suppose that $\CF$ is a weak Q-filter and we show that any
strategy \$ for \BI is not a winning strategy.  Define a sequence of
integers $\la \pi_k : k \in \omega \ra$ such that $[\pi_0,\infty) \sub
\$(\emptyset)$ and more generally

\[ [\pi_{k+1}, \infty) \sub \bigcap \{\$(\la n_0,n_1,\cdots,n_i \ra):
       n_0<n_1<\cdots<n_i<\pi_k \}. \]

\noi By assumption there is an $X \in \CF^+$ which meets each interval
$[\pi_k,\pi_{k+1})$ in at most one point. But $X=X_0 \cup X_1$ where
$X_i= X \cap \bigcup_k [\pi_{2k+i}, \pi_{2k+i+1})$ and therefore
 $X_i \in \CF^+$ for some $i$.  Write $X_i$ in increasing order as $\la
n_k : k \in \omega \ra$; but then $X_i$ is the outcome of a legal play
won by \BII, and thus \$ was not a winning strategy for \BI. 

Now we deal with player \BII.  Suppose that $\CF$ is
$\omega$-diagonalized by $\la X_n:n \in \omega \ra \sub \so$.  Fix a
surjective map $\sigma :\omega \rt \omega$ such that the preimage of
each $n$ is infinite.  Here is \BII's strategy: at stage $k$, after \BI
played $Y_k \in \Fr$, \BII responds with an element of $Y_k \cap
X_{\sigma(k)} \sm k$.  At the end of the play, \BII's outcome is a set
with infinite intersection with each $X_n$ and therefore belongs to
$\CF^+$, thus this is a winning strategy for \BII. 

\noi Now let \$ be a winning strategy for \BII in the game, we show that
$\CF$ is $\omega$-dia\-go\-na\-li\-za\-ble. We claim that

\[ (\fa Y \in \CF)(\ex n=n(Y))(\ex s=s(Y)
      \in ^n \omega)(\fa^\infty t \in ^{n+1} \omega)
       s < t \implies \$(t) \in Y. \]

Indeed otherwise one quickly produces a winning play for \BI.  But then
the collection $\{ \{ \$(s^n); n \in \omega \}; s \in \ffo \}$
$\omega$-diagonalize $\CF$. This completes the proof. \qed

\bigskip

\noi \underline{{PROOF OF THEOREM  \ref{FR-W-FC}:}} This is trivial; \BII 
chooses an infinite  $X \notin \CF$ and plays continually members of $X$. \qed

\bigskip

\noi \underline{{PROOF OF THEOREM  \ref{FR-W-F+C}:}}  If $\CF=\Fr$, then 
$\CF^+ = \so$ and \BI's strategy is to ensure that \BII's outcome is infinite,
and of course \BII has no winning strategy. 

\noi If however $\CF \neq \Fr$, \BII chooses an infinite $X \notin
\CF^+$ and continually plays members of $X$; this constitutes a winning
strategy for \BII.  \qed
 
\bigskip

\noi \underline{{PROOF OF THEOREM  \ref{DUAL-FR}:}} This follows from
Theorem \ref{BLASS} as $\Fr^+=\so$. \qed

\bigskip

\noi \underline{{PROOF OF THEOREM  \ref{F-W-F}:}} If \BII had a
winning strategy in this game, it would also be a winning strategy for
${\frak G}(\Fr, \omega, \CF)$, which is impossible.

Now for player \BI.  If $\CF$ is not a Q-filter, then \BI fixes a
partition of $\omega$ into finite sets $\la s_k : k \in \omega \ra$ and
plays to ensure that \BII's outcome meets each $s_k$ in at most one
point; therefore \BI wins.  If on the other hand $\CF$ is not a
P-filter, \BI then fixes a sequence $\la X_n : n \in \omega \ra \sub
\CF$ such that no $X \in \CF$ is almost included in each $X_n$.  It
suffices for \BI to play $\cap_{i<k}X_i$ at stage $k$ to produce a
winning strategy. 

\noi This leaves us with the more interesting situation in which $\CF$ is
a a Q-filter P-filter and we must show that any strategy \$ for \BI cannot
be a winning one. 

\noi Fixing such a strategy \$, \BII first chooses $Y \in \CF$ such that $Y
\sub^* \$(s)$ for all $s \in \ffo$.  Now \BII defines a sequence $\la
X_n : n \in \omega \ra \sub \CF$ and a sequence $\la n_k : k \in \om
\ra$ as follows: $X_0=\$(\emptyset)$ and $n_0$ is such that $Y \sm n_0
\sub X_0$. 
Now given $n_0<n_1<\cdots<n_k$, put

\[X_{k+1} = \bigcap \{ \$(\la m_0,m_1,\cdots,m_i\ra): i\leq k, 
               m_j \in [n_j,n_k] \cap Y \}. \]

\noi Then \BII chooses $n_{k+1}$ such that $Y\sm n_{k+1} \sub X_{k+1}$. 
Now because $\CF$ is a Q-filter, \BII knows very well that there is a
set in $\CF$ missing infinitely many intervals, say
$[n_{k_n},n_{{k_n}+1})$, for an infinite set $\{k_n : n \in \om \}$, and
therefore by selectivity again \BII can find a set $Y'=\{y_n : n \in
\om\} \sub Y$ in $\CF$ such that $Y'$ misses all these intervals
$[n_{k_n},n_{{k_n}+1})$ and further intersects each other interval
$[n_{k_n},n_{k_{n+1}})$ in at most one point. 

\noi But we claim now that $Y'$ is a legal play in the game! Indeed, for
each $k$, say $y_k \in [n_{k_n},n_{k_{n+1}}) \cap Y' =
[n_{{k_n}+1},n_{k_{n+1}}) \cap Y'$, then $y_k \in X_{{k_n}+1} \sub
\$(\la y_0,y_1,\cdots,y_{k-1} \ra)$.  Thus \BI's strategy \$ was
definitely not a winning one.  \qed

\bigskip

\noi \underline{{PROOF OF THEOREM  \ref{F-W-F+}:}} We first look
at player \BI. If $\CF$ is not weakly Ramsey, there is an $\CF$-tree $\CT$ such that
no branches belong to $\CF^+$; this gives a winning strategy for \BI by
playing along the tree.

\noi Now suppose that $\CF$ is weakly Ramsey and that \$ is a strategy
for \BI, we will produce a winning play for \BII showing that \BI cannot
have a winning strategy. Let $X_0 = \$(\emptyset)$ and choose $n_0 \in X_0$.
Having produced $n_0<n_1<\cdots <n_k$, let
\[X_{k+1} = \bigcap \{ \$(y_0,y_1,\cdots,y_i) :
   i \leq k, y_0 < y_1 < \cdots y_k \leq n_k \} \sm n_k. \]

\noi $X_{k+1} \in \CF$ and choose $n_{k+1} \in X_{k+1}$.  Now define a
$\CF$-tree $\CT$ inductively by letting $X_\emptyset = X_0$, and given
$s=\la y_0,y_1,\cdots , y_i \ra \in \CT$, say $n_{k-1} \leq y_i < n_k$,
then let $X_s = X_{k+1}$.  This $\CF$-tree $\CT$ must contain a branch
$\{y_k : k \in \om \} \in \CF^+$ by assumption, but this clearly is a
legal play of the game in which \BII wins. 

Now we deal with player \BII.  Suppose that $\CF$ is
$\omega$-+-diagonalized by $\la X_n:n \in \omega \ra \in \CF^+$.  Fix a
surjective map $\sigma :\omega \rt \omega$ such that the preimage of
each $n$ is infinite.  Here is \BII's strategy: at stage $k$, after \BI
played $Y_k \in \CF$, \BII responds with an element of $Y_k \cap
X_{\sigma(k)} \sm k$.  Notice the importance of having each $X_n \in
\CF^+$.  At the end of the play, \BII's outcome is a set with infinite
intersection with each $X_n$ and therefore belongs to $\CF^+$, thus it
is a winning strategy for \BII. 

\noi Now let \$ be a winning strategy for \BII in the game, we show that
$\CF$ is $\omega$-+-diagonalizable. We first define an $\CF^+$-tree as follows.
Let $X_\emptyset = \{ \$(X): X \in \CF \} \in \CF^+$, and for each $n \in 
X_\emptyset$, select an $X^n_\emptyset \in \CF$ such that
\$$(X^n_\emptyset)=n$. More generally, given $X^n_{\ov s} \in \CF$, 
for $n \in X_{\ov s} \in
\CF^+$, say $ {\ov s}=\la s_0,s_1,\cdots,s_i \ra$, 
let 
\[ X_{{\ov s} \con n} = \{ \$(X^{s_0}_\emptyset,X^{s_1}_{\la s_0 \ra},
  X^{s_2}_{\la s_0,s_1 \ra}, \cdots, X^{s_i}_{\la s_0,s_1,\cdots,s_{i-1} \ra},
         X) : X \in \CF \} \in \CF^+, \]

\noi and for each $k \in X_{{\ov s} \con n}$, choose
 $X^k_{{\ov s} \con n} \in \CF$ such that 

\[ \$(X^{s_0}_\emptyset,X^{s_1}_{\la s_0 \ra},
   X^{s_2}_{\la s_0,s_1 \ra}, \cdots, X^{s_i}_{\la s_0,s_1,\cdots,s_{i-1} \ra},
   X^k_{{\ov s} \con n})=k. \]

\noi Therefore we obtain a $\CF^+$-tree $\CT$ each of whose branches is a
legal play of the game, and therefore belongs to $\CF^+$ as \$ was a
winning strategy for \BII and this means that the sets $\{X_{\ov s}: {\ov s} \in
\CT\} \sub \CF^+$ must $\omega$-+-diagonalize $\CF$. \qed

\bigskip

\noi \underline{{PROOF OF THEOREM  \ref{F-W-FC}:}}
As in ${\frak G}(\Fr, \omega, \CF^c)$, \BI has no winning strategy.

As for \BII, let us consider the case when $\CF$ is not a Ramsey
ultrafilter first.  There are three possibilities which \BII figures
out.  If $\CF$ is not an ultrafilter, \BII chooses $X \in \CF^+ \sm \CF$
and constantly plays members of $X$, and therefore wins the game. 
Otherwise \BII checks whether $\CF$ is a Q-filter and if not chooses a
partition $\la s_k: k \in \om \ra$ of $\om$ into finite sets for which
$\CF$ contains no selector.  But then \BII's strategy is to play a
selector for the partition, therefore winning again.  Finally \BII
realizes that it must be that $\CF$ is not a P-filter and thus selects
$\la X_n :n \in \om \ra \sub \CF$ with no $Y \in \CF$ almost included in
each $X_n$.  Then \BII's strategy at the $k^{th}$ move is to play a
member of $\cap_{i<k}X_i$ and again has a winning strategy. 

\noi Now we suppose that $\CF$ is a Ramsey ultrafilter and must show
that \BII cannot have a winning strategy. But as $\CF$ is an ultrafilter,
$\CF=\CF^+$ and by the duality theorem \ref{BLASS}, the
game  ${\frak G}(\CF , \omega, \CF^c)$ is dual to the game
 ${\frak G}(\CF, \omega, \CF)$, and as $\CF$ is a Q-filter P-filter, \BI
has no winning strategy in ${\frak G}(\CF, \omega, \CF)$ by theorem
\ref{F-W-F} and therefore \BII has no winning strategy in ${\frak G}(\CF,
\omega, \CF^c)$.  \qed

\bigskip

\noi \underline{{PROOF OF THEOREM  \ref{F-W-F+C}:}} As far as player 
\BI is concerned, if $\CF$ is generated by countably many sets $\la X_n :
n \in \om \ra$, then it suffices for \BI to play 
 $\cap_{i<k}X_i \sm k$ at stage $k$, and thus the outcome of the play is a set $Y
\sub^* X_n$ for each $n$, definitely in $\CF^+$ and \BI wins.

\noi If on the other hand \BI has a winning strategy \$ in the game,
then the filter $\CF$ must be generated by $\la \$(s): s \in \ffo \ra $,
as can easily be verified.

Now for player \BII. If $\CF$ is not +-Ramsey, then there is an $\CF^+$-tree
$\CT$ none of whose branches belong to $\CF^+$; therefore \BII's strategy
is to play along the tree $\CT$.

\noi Suppose finally that $\CF$ is +-Ramsey and we show that any strategy \$ for 
\BII is not a winning strategy.

\noi We first define a $\CF^+$-tree as follows.
Let $X_\emptyset = \{ \$(X): X \in \CF \} \in \CF^+$, and for each $n \in 
X_\emptyset$, select an $X^n_\emptyset \in \CF$ such that
\$$(X^n_\emptyset)=n$. More generally, given $X^n_{\os} \in \CF$, for $n \in X_{\os} \in
\CF^+$, say $ \os=\la s_0,s_1,\cdots,s_i \ra$, 
let 
\[ X_{\os \con n} = \{ \$(X^{s_0}_\emptyset,X^{s_1}_{\la s_0 \ra},
         X^{s_2}_{\la s_0,s_1 \ra}, \cdots, X^{s_i}_{\la s_0,s_1,\cdots,s_{i-1} \ra},
         X) : X \in \CF \} \in \CF^+, \]

\noi and for each $k \in X_{\os \con n}$, choose $X^k_{\os \con n} \in \CF$ such
that 

\[\$(X^{\os_0}_\emptyset,X^{s_1}_{\la s_0 \ra},
   X^{s_2}_{\la s_0,s_1 \ra}, \cdots, X^{s_i}_{\la s_0,s_1,\cdots,s_{i-1} \ra},
   X^k_{\os \con n})=k. \]

\noi Therefore we obtain a $\CF^+$-tree $\CT$ all of whose branches are a
legal play of the game, but there is such a branch in $\CF^+$ as $\CF$ is
+-Ramsey and therefore \BII's strategy is not a winning strategy.
\qed

\bigskip

\noi \underline{{PROOF OF THEOREM  \ref{DUAL-F}:}} Follows immediately
from Theorem \ref{BLASS}. \qed

\bigskip

\noi \underline{{PROOF OF THEOREM  \ref{FR-FSW-F}:}} We must show
that the games ${\frak G}_1(\CF)$ is equivalent to the game ${\frak
G}(\Fr, \fso,\CF)$.  As \BI is trying to produce an outcome out of $\CF$
in the game ${\frak G}(\Fr, \fso,\CF)$, \BI plays without loss of
generality cofinite sets of the form $[n,\infty )$ which we identify
with $n$ to simplify notation. 

\noi Suppose first that \BI has a winning strategy \$ in the game  ${\frak G}_1(\CF)$,
we define a strategy $ \ov{\$}$ for \BI in the game ${\frak G}(\Fr, \fso,\CF)$ by

\[ \ov{\$}(s_0,s_1,\cdots,s_i) = \$(\bigcup_{j\leq i} s_j), \]

\noi where $(\bigcup_{j\leq i} s_j)$ is considered as an element of
$\fso$. This gives a winning strategy for \BI in ${\frak G}(\Fr, \fso,\CF)$.

\noi Now let \$ be a winning strategy for \BI in ${\frak G}(\Fr, \fso,\CF)$, 
we define a strategy $\ov {\$}$ for \BI in ${\frak G}_1(\CF)$ as follows:

\noi at stage $k$, after \BI played $m_0,m_1,\cdots , m_k$ and \BII
responded with $n_0,n_1,\cdots , n_k$, \BI puts $s_i = (m_{i-1},m_i]$
$\cap$ $\{n_0,n_1,\cdots , n_k\}$ for $i \leq k$ and $s_{k+1} = [m_k,
\infty)$ $\cap$ $\{n_0,n_1,\cdots , n_k\}$ and replies with
\$$(s_0,s_1,\cdots,s_{k+1})$.  This must be a winning strategy for \BI. 
Indeed let $m_0 = \$(\emptyset )$ be the first move of \BI.  Now given
$m_k$, there must be a stage at which \BII responded with an integer
bigger than \BI's play, as otherwise \BI wins the play trivially, so
call this last move of \BI by $m_{k+1}$.  Then if $s_k$ denotes \BII's
moves in the interval $[m_k,m_{k+1})$, then we get a legal play $\la
m_0,s_0,m_1,s_1,\cdots \ra$ in the game ${\frak G}(\Fr, \fso,\CF)$ which
\BI wins by assumptions; but as the outcome is the same we have shown
that $\ov{\$}$ is a winning strategy for \BI in the game ${\frak
G}_1(\CF)$.

Now we deal with player \BII. If \$ is a winning strategy for \BII in the game
${\frak G}(\Fr, \fso,\CF)$, then \BII behaves as follows in the game 
 ${\frak G}_1(\CF)$: at stage $k$, \BII imagines that \BI has played
$\la m_0, m_1, \cdots, m_k \ra$ in the game ${\frak G}(\Fr, \fso,\CF)$ 
and plays one by one the elements of \$$(\la m_0, m_1, \cdots, m_k \ra )$
without noticing \BI's moves until done, and then remembers \BI's last
move, $m_{k+1}$. At the end, the outcome is one of the game
 ${\frak G}(\Fr, \fso,\CF)$ and thus \BII wins.

\noi Finally suppose that \$ is a winning strategy for \BII in the game 
${\frak G}_1(\CF)$. 

\begin{claim} \cite{BS}
For each $\sigma \in \fso$, there is $\tau \in \fso$ and $n \in \om$,
such that \$$(\sigma \con \tau \con m) > m$ for each $m >n$.
\end{claim}

\noi If the claim were false, then it is left to the reader to verify that
\BI could easily win a play. 

\noi So here is what \BII does in the game ${\frak G}(\Fr, \fso,\CF)$:
at stage $k$, after \BI has played $\sigma = \la m_0, m_1, \cdots , m_k
\ra$, \BII chooses the corresponding $\tau = \la t_0, t_1, \cdots,
t_\ell \ra$ and $n_k$ from the lemma and replies with $s_k = \$ " \sigma
\con (\tau \rest i)$ for $i \leq \ell$; when \BI responds with
$m_{k+1}$, \BII just imagines that \BI actually responded with some
$\ov{m}_{k+1} > n_k$.  This gives a simiulated play from ${\frak
G}_1(\CF)$ and therefore \BII wins.  \qed

\bigskip

\noi \underline{{PROOF OF THEOREM  \ref{FR-FSW-F+}:}} The fact that this
game is dual to ${\frak G}(\Fr, \fso,\CF)$ is purely accidental I believe.

By playing two games simultaneously, \BII can produce two outcomes whose
union is cofinite, therefore one of them must be in $\CF^+$ and thus
\BI has no winning strategy. More precisely, let $m_0=\ov{m}_0$ be \BI's first
move as we again identify a cofinite set $[m, \infty)$ with $m$.
\BII replies with $\{ m_0 \}$, and then \BI responds with $m_1$ in the first 
game. Now  in the second game \BII replies with $[m_0,m_1]$ and waits for
\BI's response $\ov{m}_1$; then \BII comes back to the first game and replies
with $[m_1, \ov{m}_1]$. Continuing this way, \BII produces the outcome
$A = \bigcup_i [m_i, \ov{m}_i]$ in the first game, and $B=\bigcup_i [\ov{m}_i,
m_{i+1}]$ in the second. One of these sets is in $\CF^+$ and therefore \BI
lost one of the games.

Now for player \BII. If $\CF$ is meager, then \BII has definitely an easy
time winning the game; indeed there must be a sequence $\pi_0 < \pi_1 <
\cdots $ such that each member of $\CF$ meets all but finitely many of the 
intervals $[\pi_k, \pi_{k+1})$. Therefore at stage $\ell$, \BII plays one
of the intervals $[\pi_k, \pi_{k+1})$ with $k > \ell$.

\noi So we must show this is the only way that \BII can have a winning strategy.
So fix such a winning  strategy \$ for \BII. Define a sequence of integers
by $\pi_0=1$, and given $\pi_k$, choose

\[ \pi_{k+1} = max \{ \$(m_0, m_1, \cdots ,m_i ):
         i \leq k \mbox{ and } m_j \leq \pi_k \}+1. \]

\noi then each member of $\CF$ must meet all but finitely many of the intervals
$[\pi_k,\pi_{k+1})$. Otherwise, say $Y \in \CF$ misses the intervals
$[\pi_{k_\ell},\pi_{{k_\ell}+1})$ for $\ell \in \om$, then \BI wins by 
playing exactly these $k_\ell$. \qed

\bigskip

\noi \underline{{PROOF OF THEOREM  \ref{FR-FSW-FC}:}} Follows immediately
by Lemma \ref{EQUI} and Theorem \ref{FR-W-FC}. \qed

\bigskip

\noi \underline{{PROOF OF THEOREM  \ref{FR-FSW-F+C}:}} Follows immediately
by Lemma \ref{EQUI} and Theorem \ref{FR-W-F+C}. \qed

\bigskip

\noi \underline{{PROOF OF THEOREM  \ref{F-FSW-F}:}} If \BII had a
winning strategy in this game, it would be a winning strategy in the game
${\frak G}(\Fr, \fso, \CF)$, which is impossible.

If $\CF$ is meager, then \BI uses the same strategy as for 
${\frak G}(\Fr, \fso, \CF)$; if on the other hand $\CF$ is not a P-filter,
then \BI chooses a witness $\la X_n : n \in \omega \ra $ and plays 
$\bigcap_{i<k}X_i$ at stage $k$ which provides a winning strategy.

\noi Now, as a strategy for player \BI is nothing else but an
$\CF$-tree of finite sets, it should be rather clear that \BI has no
winning strategy if and only if any $\CF$-tree has a branch in $\CF$,
that is if and only if $\CF$ is a non-meager P-filter by Lemma \ref{NMP}.
\qed

\bigskip

\noi \underline{{PROOF OF THEOREM  \ref{F-FSW-F+}:}} As far as player \BI
is concerned, the proof is entirely similar to that of Theorem  \ref{F-W-F+}.

So we consider the situation for \BII.  First fix $\CF$-universal sets
$\la X_n : n \in \om \ra$ $\om$-diagonalizing $\CF$.  Player \BII fixes
a surjection $\sigma:\om \rightarrow \om$ such that the preimage of
every $n$ is infinite and at stage $k$, after \BI produced a set $Y_k
\in \CF$, \BII responds with $s_k \in X_{\pi(k)} \cap [Y_k \sm
k]^{<\om}$.  At the end of the play, \BII has produced $S= \bigcup_k
s_k$ which contains infinitely many members of each $X_n$, and is
therefore in $\CF^+$. 

\noi Now let \$ be a winning strategy for \BII and we define a tree $\CT
\sub \ffso$ such that the successors of each node ${\ov s} \in \CT$ form an
$\CF$-universal set and the collection $\la X_{\ov s}: {\ov s} \in \CT \ra$
$\om$-diagonalizes $\CF$. 

\noi Let $X_\emptyset = \{ \$(X): X \in \CF \}$, an $\CF$-universal set,
and for each $s \in X_\emptyset$, choose $X^s_\emptyset \in \CF$ such 
that \$$(X^s_\emptyset) = s$. In general, given $X^t_{\ov s} \in \CF$,
 for ${\ov s}=\la s_0,s_1, \cdots, s_i \ra \in \ffso$,
define 

\[ X_{{\ov s}\con t} = \{ \$(X^{s_0}_\emptyset, X^{s_1}_{s_0}, \cdots ,
             X^t_{\ov s}, X)): X \in \CF \} \]

\noi a $\CF$-universal set, and for each $u \in  X_{{\ov s}\con t}$, choose
$X^u_{{\ov s}\con t} \in \CF$ such that 

\[ \$(X^{s_0}_\emptyset, X^{s_1}_{s_0}, \cdots ,
             X^t_{\ov s}, X^u_{{\ov s}\con t}) = u. \]

\noi Now every branch of $\CT$ constitutes an outcome of a play of the
game, and therefore the sets ($\CF$-universal) $X_{\ov s}$ for ${\ov s}
\in \CT$ must $\om$-diagonalize $\CF$.  \qed

\bigskip

\noi \underline{{PROOF OF THEOREM  \ref{F-FSW-FC}:}} Follows immediately from
Lemma \ref{EQUI} and Theorem \ref{F-W-FC}. \qed

\bigskip

\noi \underline{{PROOF OF THEOREM \ref{F-FSW-F+C}:}} Follows immediately
from Lemma \ref{EQUI} and Theorem \ref{F-W-F+C}.  \qed

\bigskip

\noi \underline{{PROOF OF THEOREM  \ref{F+-FSW-F}:}} Player \BII
definitely has no hope for a winning strategy; if there is $X \in \CF^+
\sm \CF$, then \BII is doomed, and otherwise $\CF = \CF^+$ and the result
follows by Theorem \ref{F-FSW-F}.

So player \BI has a winning strategy if $\CF^+ \neq \CF$, and otherwise,
the result follows again from Theorem \ref{F-FSW-F}. \qed

\bigskip

\noi \underline{{PROOF OF THEOREM  \ref{F+-FSW-F+}:}}  The proof is entirely 
similar to that of Theorem \ref{F-FSW-F+} and is left to the reader.

\bigskip

\noi \underline{{PROOF OF THEOREM  \ref{F+-FSW-FC}:}} Follows immediately from
Lemma \ref{EQUI}, Theorem \ref{BLASS} and Theorem \ref{F-W-F}. \qed

\bigskip

\noi \underline{{PROOF OF THEOREM  \ref{F+-FSW-F+C}:}} Follows immediately from
Lemma \ref{EQUI}, Theorem \ref{BLASS} and Theorem \ref{F-W-F+}. \qed

\section{Conclusion}

The combinatorial properties of a filter being meager, a P-filter or a
Q-filter have probably been the most popular; another related property
which has been around for some times is that of a filter $\CF$ being
`rapid', i.e.  for any partition of $\om$ into finite sets $\la s_k : k
\in \om \ra$, $\CF$ contains a set $X$ such that $X \cap s_k$ has size
at most $k$ for each $k$.  If we modify the games above so that at stage
$k$ player \BII responds with a finite set of size at most $k$ (or
bounded by a fixed ``unbounded'' function), then the characterizations
of the winning strategies for either player will involve modifications
of this property of being rapid, such as weakly rapid and so on.  The
details are left to the interested reader.

However, we have little information on the following interesting
variation of the games:

\noi {  Problem:} Characterize winning strategies for the games in
which player \BII responds with members of $\CF^*$.

\bibliographystyle{amsplain}

\end{document}